\documentclass[12pt,a4paper,reqno]{amsart}

\usepackage{amssymb}
\usepackage{amscd}
\usepackage[pdftex,pdfpagelabels]{hyperref}
\usepackage{enumerate}
\usepackage{graphicx}
\usepackage{siunitx}
\usepackage{tikz-cd}
\usepackage{stix}
\usepackage{bm}
\DeclareMathAlphabet\mathbfcal{LS2}{stixcal}{b}{n}
\numberwithin{equation}{section}

\usepackage{mathtools}
\usepackage[tableposition=top]{caption}
\usepackage{booktabs,dcolumn}

\DeclareFontFamily{OT1}{rsfs}{}
\DeclareFontShape{OT1}{rsfs}{n}{it}{<-> rsfs10}{}
\DeclareMathAlphabet{\mathscr}{OT1}{rsfs}{n}{it}

\theoremstyle{plain}

\newtheorem{theorem}{Theorem}[section]
\newtheorem{proposition}[theorem]{Proposition}

\theoremstyle{definition}

\newtheorem{remark}[theorem]{Remark}

\renewcommand\P{\mathbb{P}}
\newcommand\E{\mathbb{E}}

\newcommand\R{\mathbb{R}}

\newcommand\N{\mathbb{N}}

\newcommand\Q{\mathbb{Q}}

\newcommand\eps{\varepsilon}

\begin{document}

\title{On product representations of squares}

\author{Terence Tao}
\address{UCLA Department of Mathematics, Los Angeles, CA 90095-1555.}
\email{tao@math.ucla.edu}



\begin{abstract}  Fix $k \geq 2$.  For any $N \geq 1$, let $F_k(N)$ denote the cardinality of the largest subset of $\{1,\dots,N\}$ that does not contain $k$ distinct elements whose product is a square.  Erd\H{o}s, S\'ark\" ozy, and S\'os showed that $F_2(N) = (\frac{6}{\pi^2}+o(1)) N$, $F_3(N) = (1-o(1))N$, $F_k(N) \asymp N/\log N$ for even $k \geq 4$, and $F_k(N) \asymp N$ for odd $k \geq 5$.  Erd\H{o}s then asked whether $F_k(N) = (1-o(1)) N$ for odd $k \geq 5$.  Using a probabilistic argument, we answer this question in the negative.
\end{abstract}

\maketitle


\section{Introduction}

\begin{figure}
\centering
\includegraphics[width=0.5\textwidth]{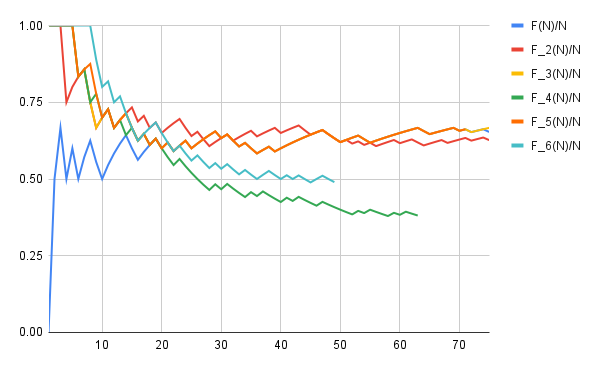}
\caption{A plot of $F_k(N)/N$ and $F(N)/N$ for $2 \leq k \leq 6$ and $N \leq 75$, though data for $k=4,5,6$ are only available up to $N=63,71, 49$ respectively.  Numerically, we have $F_3(N)=F(N)$ for $18 \leq N \leq 74$ (but not $N=75$) and $F_5(N)=F(N)$ for all $N \geq 18$ for which data is available. It is known that $F(N)/N$, $F_2(N)/N$, $F_3(N)/N$, $F_4(N)/N$, and $F_6(N)/N$ are converging to $1-c = 0.8284\dots$, $\frac{6}{\pi^2} = 0.6079\dots$, $1$, $0$, and $0$ respectively. The limiting behavior of $F_5(N)/N$ remains open, although we prove that it does not converge to $1$.  Thanks to Michael Branicky, David Corneth, Martin Ehrenstein, and Paul Muljadi for the data, which is available on the OEIS at the indicated links.}
\label{fig:chart}
\end{figure}

For any natural numbers\footnote{In this paper, the natural numbers will start from $1$.} $N, k$, let $F_k(N)$ denote the cardinality of the largest subset $A$ of $\{1,\dots,N\}$ that does not contain $k$ distinct elements whose product is a square.  Trivially we have $F_1(N) = N - \lfloor \sqrt{N} \rfloor$ (\url{https://oeis.org/A028391}), so we assume henceforth that $k \geq 2$.  The quantity $F_k(N)$ was studied by Erd\H{o}s, S\'ark\" ozy, and S\'os \cite{ess}, who established the following asymptotic bounds as $N \to \infty$:
\begin{itemize}
  \item[(i)]  One has $F_2(N) = \left(\frac{6}{\pi^2}+o(1)\right) N$. (\url{https://oeis.org/A013928})
  \item[(ii)]  One has $F_3(N) = (1-o(1)) N$. (\url{https://oeis.org/A372306})
  \item[(iii)]  If $k$ is divisible by $4$, then $F_k(N) = (1+o(1)) \frac{N}{\log N}$.  (\url{https://oeis.org/A373119} for $k=4$)
  \item[(iv)]  If $k \geq 6$ is equal to $2 \hbox{ mod } 4$, then $F_k(N) = (\frac{3}{2}+o(1)) \frac{N}{\log N}$. (\url{https://oeis.org/A372306} for $k=6$)
  \item[(v)]  If $k \geq 5$ is odd, then $(\log 2 - o(1)) N \leq F_k(N) \leq N - \delta_k \frac{N}{\log^2 N}$ for some constant $\delta_k>0$. (\url{https://oeis.org/A373178} for $k=5$)
\end{itemize}
See Figure \ref{fig:chart}.
Here as usual we use $o(1)$ to denote a quantity that goes to zero as $N \to \infty$ (holding $k$ fixed).
In fact, these authors obtained sharper bounds on $F_k(N)$; we refer the reader to \cite[Theorems 1-7]{ess} for the precise statements.  The claim $F_4(N) = o(N)$ (which is clearly implied by (iii)) was established by Erd\H{o}s more than five decades earlier, in \cite{erdos-old}.
There have been several refinements to the bounds for even $k \geq 4$, mainly using methods from graph (and hypergraph) theory; see \cite{gyori}, \cite{sarkozy}, \cite{naor}, \cite{pach1}, \cite{pach2}, \cite{pach3}.  As the case of even $k$ is not our primary concern here, we refer the reader to \cite{pach3} or \cite{fjkps} for a summary of the most recent bounds in this area.  In \cite{verstraete} the analogous question in which the squaring function $n \mapsto n^2$ is replaced by a more general polynomial is considered; very recently, new results were obtained for the cubing function $n \mapsto n^3$ in \cite{fjkps}.

Now we turn to the case of odd $k \geq 5$.  Let $F(N)$ denote the cardinality of the largest subset $A$ of $\{1,\dots,N\}$ that does not contain an odd number of elements whose product is a square (or equivalently\footnote{Indeed, if such a function equals $-1$ on $A$, then clearly no odd product of elements of $A$ can be a square.  Conversely, observe that modulo the squares $(\Q^\times)^2$ of non-zero rational numbers, the set of even products of $A$ form a subgroup of the $2$-torsion abelian group $\Q^\times / (\Q^\times)^2$, and the set of odd products form a coset of that subgroup, which is non-trivial iff they avoid squares.  Assuming the odd products do avoid squares, then Zorn's lemma yields an index two subgroup of $\Q^\times / (\Q^\times)^2$ that avoids the odd products; but any such subgroup induces a $\{-1,+1\}$-valued completely multiplicative function.}, there is a $\{-1,+1\}$-valued completely multiplicative function that equals $-1$ on $A$), then clearly $F(N) \leq F_k(N)$ for all odd $k$.  From the equivalent definition of $F(N)$, we see that $F(N)$ is the best constant such that
$$ \sum_{n \leq N} f(n) \geq 2F(N) - N$$
for all completely multiplicative functions $f \colon \N \to \{-1,+1\}$.  The sequence $N \mapsto 2F(N) - N$ begins as
$$ 	1, 0, -1, 0, -1, 0, -1, -2, -1, 0, -1, -2, -3, -4, -3, -2, -3, -4, -5, \dots$$
(\url{https://oeis.org/A360659}).  The sequence $F(N)$ is thus
$$ 0, 1, 2, 2, 3, 3, 4, 5, 5, 5, 6, 7, 8, 9, 9, 9, 10, 11, 12, 12, \dots$$
(\url{https://oeis.org/A373114}).  The results in \cite[Theorem 1, (1.1), (1.2)]{gs}, when phrased in our notation, assert the asymptotic
\begin{equation}\label{fn}
   F(N) = (1-c+o(1)) N
\end{equation}
as $N \to \infty$, where $c$ is the Hall--Montgomery constant
\begin{equation}\label{17}
  c \coloneqq 1 - \log(1+\sqrt{e}) + 2 \int_1^{\sqrt{e}} \frac{\log t}{t+1}\ dt = 0.171500\dots
\end{equation}
(\url{https://oeis.org/A143301}).
This answered a conjecture of Hall \cite[\S 2]{hall} (and independently, Hugh Montgomery), and also answers a question in \cite{erdos}; as noted in \cite{erdos}, the weaker claim $F(N) \leq (1-c_0+o(1))N$ for some constant $c_0>0$ was observed in \cite[Theorem 4.1]{ruzsa} to follow\footnote{A short proof is as follows (Imre Ruzsa, private communication).  If $\sum_{p \leq N} \frac{1-f(p)}{p}$ is bounded away from zero then the claim follows from the Hal\'asz inequality \cite{halasz}, while if $\sum_{p \leq N} \frac{1-f(p)}{p}$ is small then a simple union bound will adequately control the values where $f(n)=-1$.} by adapting the proof of \cite[Theorem 1]{ers}; it was also showed independently in \cite{hall}, and had previously been conjectured by Heath-Brown.  In view of \eqref{fn}, the constant $\log 2 = 0.6931\dots$ in claim (v) above can be improved to $1-c = 0.828499\dots$; this was also recently observed in \cite[Theorem 9]{fjkps}.  We remark that an explicit example\footnote{Indeed, from \cite[Corollary 1]{gs} we see that this set $A$ is ``essentially'' the only example, as any other set $A'$ with $F(N)+o(N)$ elements that avoids squares in its odd products would then have to lie in the set $\{n:f'(n)=-1\}$ for some completely multiplicative function $f'$ with $\sum_{p < N^{1/(1+\sqrt{e})}} \frac{1-f'(p)}{p} + \sum_{N^{1/(1+\sqrt{e})} \leq p \leq N} \frac{1+f'(p)}{p} = o(1)$.  This hints at a possible theory of ``pretentious sets'', analogous to the theory of ``pretentious multiplicative functions'' \cite{pretentious}.  We thank Andrew Granville for this remark.} of a set $A$ that achieves the asymptotic $(1-c+o(1))N$ can be described as the set of numbers in $\{1,\dots,N\}$ with exactly one prime factor larger than $N^{1/u}$, where $u$ is optimized as $u \coloneqq 1 + \sqrt{e}=2.6487\dots$; see \cite[\S 2]{hall}, \cite[(1.2)]{gs} (as well as the remark in \cite[p. 585]{ess}) for further discussion. As noted in \cite{ess}, the choice $u=2$ would instead give the weaker lower bound $F(N) \geq (\log 2 - o(1)) N$, which explains the lower bound in claim (v) above.  The even weaker bound $F(N) \geq (\frac{1}{2}-o(1)) N$ is also immediate from the prime number theorem $\sum_{n \leq N} \lambda(n)=o(N)$ after setting $f$ to be the Liouville function $\lambda$.

\begin{remark}  One can also study the logarithmic analogue $L_k(N)$ of $F_k(N)$, defined as the largest value of $\sum_{n \in A} \frac{1}{n}$ where $A$ is a subset of $\{1,\dots,N\}$ such that no $k$ distinct elements of $A$ multiply to a square.  These quantities are better understood than $F_k(N)$ when $k$ is odd. Indeed, it was shown in \cite[Theorem 8]{ess} that $L_k(N) = (1 + o(1)) \log\log N$ when $k$ is a multiple of $4$, $L_k(N) = (\frac{3}{2}+o(1)) \log\log N$ when $k \geq 4$ is equal to $2 \hbox{ mod } 4$, and that $L_k(N) = (\frac{1}{2}+o(1)) \log N$ when $k$ is odd, say for $N \geq 10$.  The same arguments also show that the logarithmic analogue $L(N)$ of $F(N)$ is equal to $(\frac{1}{2}+o(1)) \log N$ (the bound is attained when $A = \{n \leq N: \lambda(n)=-1\}$ consists of numbers up to $N$ that are products of an odd number of primes).  It is also easy to establish that $L_2(N) = (\frac{6}{\pi^2}+o(1)) \log N$, since modulo the square numbers, $A$ can represent each squarefree number at most once.
\end{remark}

For each $k$, let $c_k^-, c_k^+$ denote the best constants for which one has the bounds $(1 - c^+_k - o(1)) N \leq F_k(N) \leq (1 - c^-_k+o(1)) N$ as $N \to \infty$, that is to say
$$ c^-_k = \liminf_{N \to \infty} 1 - \frac{F_k(N)}{N}; \quad c^+_k = \limsup_{N \to \infty} 1 - \frac{F_k(N)}{N}.$$
From \eqref{fn} and the trivial bound $F_k(N) \leq N$ we see that $0 \leq c^-_k \leq c_k^+ \leq c$ for any odd $k \geq 3$, while from claims (i)-(iv) we have $c^-_2 = c^+_2 = 1-\frac{6}{\pi^2} = 0.39207\dots$, $c^-_k = c^+_k = 0$ for $k=1,3$, and $c^-_k = c^+_k = 1$ for any even $k \geq 4$.  In \cite[Problem 1]{erdos} it was asked\footnote{\url{https://www.erdosproblems.com/121}} whether $c^-_k=c^+_k=0$ for all odd $k \geq 5$, or equivalently if $F_k(N) = (1-o(1)) N$ as $N \to \infty$ for all such $k$.  In this note, we answer this question in the negative:

\begin{theorem}[Main theorem]\label{main-thm}  We have $c^+_k \geq c^-_k>0$ for all $k \geq 4$. In other words, for any $k \geq 4$, there exists a positive quantity $c_k$ such that $F_k(N) \leq (1-c_k+o(1)) N$ as $N \to \infty$.
\end{theorem}

Interestingly, in contrast to previous results, the parity of $k$ plays no role in our arguments, although our result is only new in the case of odd $k$.  We also do not use any tools from graph theory; indeed, no combinatorial tools more advanced than double counting are used (which we found convenient to phrase in a continuous setting, via the Fubini--Tonelli theorem). Similarly, we do not need any tool from analytic number theory more advanced than the Mertens theorems and the prime number theorem; in particular, we will not need to analyze Wirsing type integral equations, as was done in \cite{gs}.

\begin{remark}\label{mono}  Our methods would in principle give an explicit lower bound for $c^-_k$ for each $k$, but it is quite small; in particular, it is much smaller than the constant $c$ defined in \eqref{17}, and deteriorates as $k \to \infty$.  However, one can improve the bounds for large $k$ by the following argument (communicated to the author by Csaba S\'andor).  One has an easy inequality
  $$ F_{k+l}(N) \leq \max(F_k(N), F_l(N)+k)$$
for any $k,l \geq 1$.  Indeed, if $A \subset \{1,\dots,N\}$ has cardinality greater than $\max(F_k(N), F_l(N)+k)$, then it contains $k$ distinct elements that multiply to a square; removing those elements, one can then locate $l$ additional distinct elements that multiply to a square, and so by combining the two we obtain $k+l$ distinct elements that multiply to a square.  Dividing by $N$ and taking limits, we conclude that
$$ c^\pm_{k+l} \geq \min( c^\pm_k, c^-_l)$$
for any $k,l \geq 1$ and choice of sign $\pm$.  In particular, since $c^\pm_k \leq c < 1-\frac{\pi^2}{6} = c^-_2$ for odd $k$, we see that $c^\pm_{k+2} \geq c^\pm_k$ for odd $k$, thus $c^\pm_k$ is monotone nondecreasing for odd $k$.  In view of the results in \cite{gs}, it seems plausible though to conjecture that $c^-_k=c^+_k$ converges to $c$ in the limit as $k$ is odd and goes to infinity; one could even more boldly conjecture that $c^-_k=c^+_k=c$ for all odd $k \geq 5$.  In view of monotonicity, it would suffice (assuming $c^-_k = c^+_k$, i.e., that $F_k(N)/N$ has a limiting value) to establish this bolder conjecture for $k=5$.  It is somewhat difficult to compute $F_k(N)$ numerically for large $N$, but the limited numerics available are not inconsistent with this conjecture; see Figure \ref{fig:chart}.
\end{remark}

The behavior of $F_k(N)$ is different in the cases  $k < 4$ and $k \geq 4$, basically because one has $\binom{k}{2} \leq k$ in the former case but not the latter.  To explain this further, we first review the case $k=3$, adapting the analysis from \cite{ess}.  Suppose $n_1,n_2,n_3$ are elements of $\{1,\dots,N\}$ such that $n_1 n_2 n_3$ is a square.  To simplify the discussion let us suppose that $n_1,n_2,n_3$ are comparable in magnitude to $N$ and are also squarefree.  Then from the fundamental theorem of arithmetic, we can factor
\begin{equation}\label{3fac}  n_1 = d_{12} d_{13}; \quad n_2 = d_{12} d_{23}; \quad n_3 = d_{13} d_{23}
\end{equation}
for some natural numbers $d_{12}, d_{13}, d_{23}$ (which will also be coprime, although this will not be relevant for our discussion).  Taking logarithms, we have
$$ \log d_{12} + \log d_{13}, \log d_{12} + \log d_{23}, \log d_{13} + \log d_{23} = \log N + O(1).$$
This is a system of $k=3$ approximate equations in $\binom{k}{2}=3$ unknowns, and is solved as
$$ \log d_{12}, \log d_{13}, \log d_{23} = \frac{1}{2} \log N + O(1).$$
In particular, $n_1,n_2,n_3$ all have a factor close to $\sqrt{N}$ in magnitude. However, from standard results\footnote{One can also invoke the Hardy--Ramanujan law here, which asserts that typical numbers of size either $\sqrt{N}$ or $N$ have $(1+o(1))\log\log N$ prime factors.} on the ``multiplication table problem'' (see, e.g., \cite{ford}), such numbers are sparse inside $\{1,\dots,N\}$, and so one expects to obtain the result (ii) above by deleting those numbers from $\{1,\dots,N\}$ (and also modifying the above argument to handle non-squarefree numbers).  We refer the reader to \cite{ess} for a more detailed treatment of this argument.

Now we discuss the $k=4$ case, in a fashion that will generalize readily to other $k \geq 4$.  Suppose that $n_1,n_2,n_3,n_4$ are elements of $\{1,\dots,N\}$ such that $n_1 n_2 n_3 n_4$ is a square.  In analogy with \eqref{3fac}, we may expect a model case to be when all the $n_i$ are comparable in size to $N$ and one has the factorization
\begin{equation}\label{4fac}  n_1 = d_{12} d_{13} d_{14}; \quad n_2 = d_{12} d_{23} d_{24}; \quad n_3 = d_{13} d_{23} d_{34}; \quad n_4 = d_{14} d_{24} d_{34}
\end{equation}
for some natural numbers $d_{12}, d_{13}, d_{14}, d_{23}, d_{24}, d_{34}$.  Taking logarithms, we now have $k=4$ approximate equations in $\binom{k}{2}=6$ unknowns, and the system is now underdetermined; in particular, there are no strong constraints on the magnitudes of the individual $d_{ij}$, or on any of the other factors of the $n_i$.  As it turns out, if one selects the $d_{ij}$ from an appropriate probability distribution, then this construction gives a random tuple $(\mathbf{n}_1,\mathbf{n}_2,\mathbf{n}_3,\mathbf{n}_4)$ multiplying to a square, which will lie in $\{1,\dots,N\}^4$ quite often, and which will be somewhat evenly distributed in the sense that no element of $\{1,\dots,N\}$ is attained with exceptionally high probability by any of the $\mathbf{n}_i$; see Proposition \ref{prob-prop} for a precise statement.  In particular, a subset of $\{1,\dots,N\}$ of density $1-o(1)$ will be quite likely to be attained by all of the $\mathbf{n}_i$, allowing one to locate elements $n_1,n_2,n_3,n_4$ of the subset that multiply to a square; it will also not be difficult to ensure that these elements are distinct.  As it turns out, this argument can be made rigorous for all $k \geq 4$, as we shall see in the next section.

\subsection{Acknowledgments}

The author is supported by NSF grant DMS-1764034. We thank Thomas Bloom, Andrew Granville,  P\'eter P\'al Pach, and Imre Ruzsa for valuable comments and help with the references.  We also thank Thomas Bloom for founding the site \url{https://www.erdosproblems.com} where the author discovered this problem, and for supplying a copy of the reference \cite{erdos}, P\'alv\"olgyi D\"om\"ot\"or for providing the recent reference \cite{fjkps}, and Will Sawin and Csaba S\'andor for a conversation leading to Remark \ref{mono}.  We also thank the anonymous referee for comments and corrections.

\section{Proof of theorem}

To prove Theorem \ref{main-thm}, we use the probabilistic method.  Fix $k \geq 4$; we use the asymptotic notation $X \ll Y$, $Y \gg X$, or $X = O(Y)$ to denote a bound of the form $|X| \leq C_k Y$ where $C_k$ depends only on $k$, and $X \asymp Y$ as an abbreviation for $X \ll Y \ll X$.  We use boldface symbols such as $\mathbf{n}$ to denote random variables, $\P$ to denote probability, and $\E$ to denote expectation.  The key construction is then

\begin{proposition}[Probabilistic construction]\label{prob-prop}  Let $N$ be sufficiently large. Then one can find a random tuple $(\mathbf{n}_1,\dots,\mathbf{n}_k)$ of natural numbers and an event $E$ obeying the following properties:
\begin{itemize}
\item[(i)] With probability $1$, $\mathbf{n}_1 \dots \mathbf{n}_k$ is a perfect square.
\item[(ii)] On the event $E$, one has $\mathbf{n}_i \leq N$ for all $i=1,\dots,k$.
\item[(iii)]  One has $\P(E) \gg \frac{1}{\log^k N}$.
\item[(iv)]  For every $1 \leq i < j \leq k$, one has $\P( \mathbf{n}_i = \mathbf{n}_j ) = o\left( \frac{1}{\log^k N} \right)$ as $N \to \infty$.
\item[(v)] For any $1 \leq n \leq N$ and $i=1,\dots,k$, one has $\P( \mathbf{n}_i = n \wedge E ) \ll \frac{1}{N \log^k N}$.
\end{itemize}
\end{proposition}

Let us now see how Theorem \ref{main-thm} follows from Proposition \ref{prob-prop}.  Suppose for contradiction that Theorem \ref{main-thm} failed for some $k \geq 4$, then one can find a sequence of $N$ going to infinity and subsets $A = A_N$ of $\{1,\dots,N\}$ with $|A| = (1-o(1)) N$ such that no product of $k$ distinct elements of $A$ is a perfect square.  Then only $o(N)$ of the elements of $\{1,\dots,N\}$ lie outside of $A$.  By Proposition \ref{prob-prop}(ii), (v), and the union bound we conclude that for each $i=1,\dots,k$, the probability that $E$ holds and $\mathbf{n}_i \not \in A$ is at most $o(\frac{1}{\log^k N})$.  From this, Proposition \ref{prob-prop} (i), (iii), (iv), and the union bound, we conclude that if $N$ is large enough in this sequence, then with positive probability, the $\mathbf{n_1},\dots,\mathbf{n_k}$ are all distinct, lie in $A$, and multiply to a square, giving the desired contradiction.

It remains to prove the proposition.  Let $\eps>0$ be a small constant depending on $k$ to be chosen later.  We subscript asymptotic notation by $\eps$ if we permit implied constants to depend on $\eps$, thus for instance $X \ll_\eps Y$ denotes the bound $|X| \leq C_{k,\eps} Y$ for some $C_{k,\eps}$ depending on both $k$ and $\eps$. Note that to prove the properties in Proposition \ref{prob-prop}, it is permissible to allow the implied constants to depend on $\eps$, since $\eps$ will ultimately depend only on $k$.

We define the following preliminary random variables:
\begin{itemize}
\item $\mathbf{d}$ is a random variable supported on the (finite) set ${\mathcal D}$ of squarefree natural numbers with all prime factors less than $N^{\eps^2}$, with probability distribution
$$ \P( \mathbf{d} = d ) \coloneqq \frac{1}{\sum_{d' \in {\mathcal D}} \frac{1}{(k-1)^{\omega(d')} d'}} \frac{1}{(k-1)^{\omega(d)} d}$$
for each $d \in {\mathcal D}$, where\footnote{All sums and products over $p$ are implicitly understood to be restricted to primes.} $\omega(d) \coloneqq \sum_{p|d} 1$ denotes the number of prime factors of $d$.
\item $\mathbf{p}$ is a random variable supported on the set ${\mathcal P}$ of primes between $N^{\eps}$ and $N$, with probability distribution
$$ \P( \mathbf{p} = p ) \coloneqq \frac{1}{\sum_{p' \in {\mathcal P}} \frac{1}{p'}} \frac{1}{p}$$
for each $p \in {\mathcal P}$.
\end{itemize}

\begin{remark} The random variable $\mathbf{p}$ has a clear interpretation; it is choosing a random prime of size very roughly of order $N$ (in doubly logarithmic scale).  To interpret the random variable $\mathbf{d}$, note that we have the Euler product
  \begin{equation}\label{euler}
  \sum_{d' \in {\mathcal D}} \frac{1}{(k-1)^{\omega(d')} d'} = \prod_{p<N^{\eps^2}} \left( 1 + \frac{1}{(k-1)p} \right).
  \end{equation}
  If we restrict attention to elements $d'$ of ${\mathcal D}$ divisible by a single prime $p$, then the factor of $1 + \frac{1}{(k-1)p}$ on the right-hand side drops to $\frac{1}{(k-1)p}$.  From these sorts of calculations one can check that for each prime $p < N^{\eps^2}$, the event that $p$ divides $\mathbf{d}$ occurs with probability $\frac{1}{(k-1)p+1}$, and that these events are independent in $p$.  Thus, while a typical number of size $N$ would be expected to have $(1+o(1))\log\log N$ prime factors by the Hardy--Ramanujan law, one can calculate that $\mathbf{d}$ should typically only have $\frac{1+o(1)}{k-1} \log\log N$ prime factors.  Standard heuristics from the anatomy of integers (in particular, the Poisson--Dirichlet distribution) then suggests that a reasonable proportion of numbers of size up to $N$ resemble the product of $k-1$ copies of $\mathbf{d}$ and $k-1$ copies of $\mathbf{p}$; this may be a useful heuristic to bear in mind in the discussion that follows.  A related heuristic that the reader may find helpful is to think of $\mathbf{d}$ as roughly similar to the number one would obtain if one started with a random (squarefree) number very roughly of order $N^{\eps^2}$ (in doubly logarithmic scale), factored it arbitrarily into the product of $k-1$ factors, and then selected one of those factors at random.
\end{remark}

From Mertens' theorem and \eqref{euler} we have
\begin{align}
  \sum_{d' \in {\mathcal D}} \frac{1}{(k-1)^{\omega(d')} d'} 
  &\asymp \log^{1/(k-1)} N^{\eps^2}\label{ab}\\
  &\asymp_\eps \log^{1/(k-1)} N\nonumber
\end{align}
and
$$ \sum_{p' \in {\mathcal P}} \frac{1}{p'} \asymp_\eps 1.$$
We conclude that
\begin{equation}\label{d-dist}
  \P( \mathbf{d} = d ) \asymp_\eps \frac{1}{\log^{1/(k-1)} N} \frac{1}{(k-1)^{\omega(d)} d}
\end{equation}
for $d \in {\mathcal D}$, and
\begin{equation}\label{p-dist}
  \P( \mathbf{p} = p ) \asymp_\eps \frac{1}{p}
\end{equation}
for $p \in {\mathcal P}$.

Now let $\mathbf{d}_{i,j}$, $\mathbf{p}_{i,j}$ for $1 \leq i < j \leq k$ be copies of $\mathbf{d}$, $\mathbf{p}$, with all of the random variables $\mathbf{d}_{i,j}$, $\mathbf{p}_{i,j}$ jointly independent.  We also define $\mathbf{d}_{j,i} \coloneqq \mathbf{d}_{i,j}$ and $\mathbf{p}_{j,i} \coloneqq \mathbf{p}_{i,j}$ for $1 \leq i < j \leq k$.  We then define the random variables $\mathbf{n}_i$ for $1 \leq i \leq k$ by the formula
\begin{equation}\label{ni-def}
  \mathbf{n}_i \coloneqq \prod_{1 \leq j \leq k: j \neq i} \mathbf{d}_{i,j} \mathbf{p}_{i,j}.
\end{equation}
Then we have
$$ \prod_{i=1}^k \mathbf{n}_i = \left( \prod_{1 \leq i < j \leq n} \mathbf{d}_{i,j} \mathbf{p}_{i,j} \right)^2$$
so property (i) of Proposition \ref{prob-prop} is satisfied.  We can also easily establish property (iv) as follows. If we had $\mathbf{n}_i = \mathbf{n}_j$ for some $1 \leq i < j \leq k$, then by dividing by $p_{i,j}$ and then comparing the largest remaining prime factors, we see that $\mathbf{p}_{i,l} = \mathbf{p}_{j,m}$ for some $l,m$ distinct from both $i$ and $j$.  But by \eqref{p-dist}
and the union bound we see that such an event occurs with probability at most $O_\eps(1/N^\eps)$, giving Proposition \ref{prob-prop}(iv) (with plenty of room to spare).

Now we let $E$ be the event that $\prod_{1 \leq i < j \leq k} \mathbf{d}_{i,j}$ is squarefree and $N/2 \leq \mathbf{n}_i \leq N$ for all $i=1,\dots,k$.  Clearly Proposition \ref{prob-prop}(ii) is obeyed.  Now we establish (iii).  We first use Rankin's trick to get some crude probabilistic upper bounds on the $\mathbf{d}_{i,j}$. From Mertens' theorem we have
\begin{align*}
\sum_{d < N^{\eps^2}} \frac{1}{(k-1)^{\omega(d)} d^{1-\frac{1}{\log N^{\eps^2}}}}
&= \prod_{p < N^{\eps^2}} \left( 1 + \frac{p^{1/\log N^{\eps^2}}}{(k-1)p} \right) \\
&\leq \exp\left( \sum_{p < N^{\eps^2}} \frac{1+O(\log p / \log N^{\eps^2})}{(k-1) p}\right) \\
&\ll \log^{\frac{1}{k-1}} N^{\eps^2}
\end{align*}
and thus by \eqref{ab}
$$ \E \mathbf{d}^{\frac{1}{\log N^{\eps^2}}} \ll 1.$$
In particular, by Markov's inequality we have that
$$ \P( \mathbf{d} \leq N^{\eps}) = 1 - O( e^{-1/\eps}).$$
By independence, we thus see that with probability $1 - O(e^{-1/\eps})$, we have
$$ \mathbf{d}_{i,j} \leq N^{\eps}$$
for all $1 \leq i < j \leq k$.  Also, for any prime $1 \leq p < N^{\eps^2}$, the probability that $\prod_{1 \leq i < j \leq k} \mathbf{d}_{i,j}$ is divisible by $p^2$ is bounded away from zero, and is $1 - O(1/p^2)$ for large $p$.  As these events are independent, we see that the probability that $\prod_{1 \leq i < j \leq k} \mathbf{d}_{i,j}$ is squarefree is $\gg 1$, if $\eps$ is small enough.

Let us condition the $\mathbf{d}_{i,j}$ to be fixed quantities $d_{i,j}$ that are bounded by $N^\eps$ and whose product is squarefree; by the above discussion, this occurs with probability $\gg 1$. Note that the $\mathbf{p}_{i,j}$ remain independent copies of $\mathbf{p}$ after this conditioning.  If we introduce the target magnitudes
$$ N_i \coloneqq N / \prod_{1 \leq j \leq k: j \neq i} d_{i,j}$$
then we have
\begin{equation}\label{nie}
 N_i = N^{1-O(\eps)}
\end{equation}
for all $i=1,\dots,k$, and by the law of total probability it will suffice to establish the bound
\begin{equation}\label{pinn}
  \P\left( N_i/2 \leq \prod_{1 \leq j \leq k: j \neq i} \mathbf{p}_{i,j} \leq N_i\ \forall i=1,\dots,k\right) \gg_\eps \log^{-k} N.
\end{equation}
We establish this by a double counting argument (in the ``continuous'' form of the Fubini--Tonelli theorem).  Let $U \subset \R^{\binom{k}{2}}$ be the polytope of all tuples $(u_{i,j})_{1 \leq i < j \leq k}$ of real numbers $u_{i,j}$ such that
\begin{equation}\label{2en}
   2\eps \log N \leq u_{i,j} \leq 0.9 \log N
\end{equation}
for all $1 \leq i < j \leq k$, and such that
\begin{equation}\label{lne}
   \log N_i - 2\eps \leq \sum_{1 \leq j \leq k: j \neq i} u_{i,j} \leq \log N_i - \eps
\end{equation}
for all $1 \leq i \leq k$, where we adopt the convention $u_{j,i} = u_{i,j}$ for $1 \leq i < j \leq k$.  Note that the collection of linear forms
\begin{equation}\label{u-form}
  (u_{i,j})_{1 \leq i < j \leq k} \mapsto \sum_{1 \leq j \leq k: j \neq i} u_{i,j}
\end{equation}
is linearly independent.  Indeed, if this were not the case, then there would exist real coefficients $a_1,\dots,a_k$, not all zero, such that the linear form
$$ \sum_{i=1}^k a_i \sum_{1 \leq j \leq k: j \neq i} u_{i,j}$$
vanished identically, which is equivalent to $a_i + a_j = 0$ for all $1 \leq i < j \leq k$.  But this easily implies (since $k \geq 3$) that all the $a_i$ vanish, a contradiction.

Because of this linear independence, the affine space of $(u_{i,j})_{1 \leq i < j \leq k}$ for which $\sum_{1 \leq j \leq k: j \neq i} u_{i,j} = \log N_i$ for all $1 \leq i \leq k$ has codimension $k$, and by \eqref{nie} it lies within $O(\eps \log N)$ of the point
$(\frac{1}{k-1} \log N)_{1 \leq i < j \leq k}$, which is well inside the interior of the region \eqref{2en}.  From this we conclude that the region $U$ has volume $\asymp_\eps \log^{\binom{k}{2}-k} N$.

The significance of this region $U$ is as follows.  If $(u_{i,j})_{1 \leq i < j \leq k} \in U$, and
\begin{equation}\label{lpu}
  |\log \mathbf{p}_{i,j} - u_{i,j}| \leq \frac{\eps}{k}
\end{equation}
for all $1 \leq i < j \leq k$ then from \eqref{lne} and the triangle inequality we have $N_i/2 \leq \prod_{1 \leq j \leq k: j \neq i} \mathbf{p}_{i,j} \leq N_i$ for all $i$.  On the other hand, from the prime number theorem and \eqref{2en}, \eqref{p-dist}, we see that for each such $(u_{i,j})_{1 \leq i < j \leq k}$, the probability that \eqref{lpu} holds is $\gg_\eps \log^{-\binom{k}{2}} N$.  Integrating this over $U$ and using the Fubini--Tonelli theorem, we conclude that the random variable
$$ \mathrm{vol} \{ (u_{i,j})_{1 \leq i < j \leq k} \in U: \eqref{lpu} \hbox{ holds} \}$$
has expectation $\gg_\eps \log^{-k} N$.  On the other hand, because the region whose volume is computed has diameter $O(1)$, this random variable is bounded by $O(1)$, hence the probability that this random variable is non-zero is $\gg_\eps \log^{-k} N$, which gives the required bound \eqref{pinn}.

The only remaining task is to establish Proposition \ref{prob-prop}(v).  By symmetry we may take $i=1$.  Fix $n=1,\dots,N$; it will suffice to show that
$$ \P( \mathbf{n}_1 = n \wedge E ) \ll_\eps \frac{1}{N \log^k N}.$$
From \eqref{ni-def} and the construction of $E$, we see that the event here is only non-empty if
$$ N/2 \leq n \leq N,$$
and $n = dq$ where $d$ is squarefree with all prime factors less than $N^{\eps^2}$, and $q$ is the product of $k-1$ primes between $N^{\eps}$ and $N$.  From \eqref{ni-def} this forces
$$\mathbf{d}_{1,2} \dots \mathbf{d}_{1,k} = d$$
and
$$\mathbf{p}_{1,2} \dots \mathbf{p}_{1,k} = q,$$
so there are $(k-1)^{\omega(d)}$ choices for the $\mathbf{d}_{1,2},\dots,\mathbf{d}_{1,k}$, and at most $O(1)$ choices for the $\mathbf{p}_{1,2}, \dots, \mathbf{p}_{1,k}$.  By \eqref{d-dist}, \eqref{p-dist}, each such choice $d_{1,2},\dots,d_{1,k}, p_{1,2},\dots,p_{1,k}$ of these variables occurs with probability
$$ \asymp_\eps \prod_{j=2}^k \frac{1}{\log^{1/(k-1)} N} \frac{1}{(k-1)^{\omega(d_{1,j})} d_{i,j}} \frac{1}{p_{i,j}} = \frac{1}{(k-1)^{\omega(d)} d q \log N} \asymp \frac{1}{(k-1)^{\omega(d)} N \log N}$$
and so the sum of all these probabilities is $\asymp_\eps \frac{1}{N \log N}$.  Thus, for each such choice, if we condition to the event that $\mathbf{d}_{1,j} = d_{1,j}$ and $\mathbf{p}_{1,j} = p_{1,j}$, it suffices to show that the event $E$ holds with conditional probability $\ll_\eps 1 / \log^{k-1} N$.

We now also condition the $\mathbf{d}_{l,j} = d_{l,j}$ to be fixed for $2 \leq l < j \leq k$.  Similarly to before, we introduce the target magnitudes
$$ \tilde N_l \coloneqq \frac{N}{p_{1,l} \prod_{1 \leq j \leq k: j \neq l} d_{l,j}}$$
for $2 \leq l \leq k$, with the convention that $d_{l,j} = d_{j,l}$ if $l > j$.  (We permit the $\tilde N_l$ to be less than $1$, as the desired bounds will be trivially true in this case.)  In view of \eqref{ni-def}, it will suffice to show that the event that
\begin{equation}\label{tl2}
   \tilde N_l/2 \leq \prod_{2 \leq j \leq k: j \neq l} \mathbf{p}_{l,j} \leq \tilde N_l\ \forall 2 \leq l \leq k
\end{equation}
holds with probability $\ll_\eps 1 / \log^{k-1} N$.

To establish this, we again use a double counting argument.  Let $\tilde U$ denote the (possibly empty) polytope of tuples  $(\tilde u_{l,j})_{2 \leq l < j \leq k}$ such that
\begin{equation}\label{2en-alt}
  \frac{\eps}{2} \log N \leq \tilde u_{l,j} \leq 2 \log N
\end{equation}
for all $2 \leq l < j \leq k$, and such that
\begin{equation}\label{lne-alt}
  \left|\sum_{2 \leq j \leq k: j \neq l} \tilde u_{l,j} - \log \tilde N_l\right| \leq 10k
\end{equation}
for all $2 \leq j \leq k$, with the convention that $\tilde u_{l,j} = \tilde u_{j,l}$ for $l > j$.  Now we make crucial use of the hypothesis $k \geq 4$ to observe that the linear forms
$$ (\tilde u_{l,j})_{2 \leq l < j \leq k} \mapsto \sum_{2 \leq j \leq k: j \neq l} \tilde u_{l,j}$$
are linearly independent; this is the same argument used to establish the linear independence of \eqref{u-form}, but with $k$ replaced by $k-1$.  In particular, the constraint \eqref{lne-alt} places $(\tilde u_{l,j})_{2 \leq l < j \leq k}$ within $O(1)$ of a codimension $k-1$ space.  From \eqref{2en-alt}, $\tilde U$ has diameter $O(\log N)$, thus $\tilde U$ has volume at most $O( \log^{\binom{k-1}{2} - (k-1)} N)$.

For any $(\tilde u_{l,j})_{2 \leq l < j \leq k}$ in $\tilde U$, we see from \eqref{p-dist} and Mertens' theorem (or the prime number theorem) that the event
\begin{equation}\label{lpu-alt}
|\log \mathbf{p}_{l,j} - \tilde u_{l,j}| \leq 1 \ \forall 2 \leq l < j \leq k
\end{equation}
holds with probability $O_\eps(\log^{-\binom{k-1}{2}} N)$.
Integrating this on $\tilde U$ and using the Fubini--Tonelli theorem, we conclude that the random variable
$$ \mathrm{vol} \{ (\tilde u_{l,j})_{2 \leq l < j \leq k} \in \tilde U: \eqref{lpu-alt} \hbox{ holds} \}$$
has expectation $O_\eps( 1 / \log^{k-1} N )$.  On the other hand, if we are in the event \eqref{tl2}, then on taking logarithms and applying the triangle inequality we see that any tuple $(\tilde u_{l,j})_{2 \leq l < j \leq k}$ obeying \eqref{lpu-alt} will lie in $\tilde U$, and so the above volume is $\gg 1$.  The claim follows.

\begin{remark}\label{many-tuples} (This remark was suggested by Andrew Granville.)  The above argument shows that if $k \geq 4$, $N$ is sufficiently large, and $A \subset \{1,\dots,N\}$ contains at least $(1 - c_k)N$ elements for some sufficiently small $c_k>0$, then $A$ contains a tuple $(n_1,\dots,n_k)$ of $k$ distinct natural numbers multiplying to a square.  In fact, a modification of the argument shows that $A$ in fact contains at least $N^{k/2} / \log^{k \log(k-1) + o(1)} N$ such tuples.  We sketch the details as follows.  We modify the event $E$ slightly by imposing the additional condition
\begin{equation}\label{om-d}
  \omega(\mathbf{d}_{i,j}) = \left(\frac{1}{k-1}+o(1)\right) \log\log N
\end{equation}
for all $1 \leq i < j \leq k$, for a suitable choice of decay rate $o(1)$; one can check from Bernstein's inequality that this does not significantly impact the arguments (and in particular one continues to have the lower bound (ii)).  Given any tuple $(d_{i,j}, p_{i,j})_{1 \leq i < j \leq k}$ with
\begin{equation}\label{dpi}
  N/2 \leq \prod_{j \neq i} d_{i,j} p_{i,j} \leq N \quad \forall 1 \leq i \leq k,
\end{equation}
we see from \eqref{d-dist}, \eqref{p-dist} that the probability that the random tuple $(\mathbf{d}_{i,j}, \mathbf{p}_{i,j})_{1 \leq i < j \leq k}$ attains this value is
$$ \ll_\eps \frac{1}{\log^k N} \prod_{1 \leq i < j \leq k} \frac{1}{(k-1)^{\omega(d_{i,j})} d_{i,j} p_{i,j}}.$$
By multiplying \eqref{dpi} over all $k$ values of $i$ and taking square roots, we have
$$ \prod_{1 \leq i < j \leq k} d_{i,j} p_{i,j} \asymp N^{k/2}$$
and hence by \eqref{om-d} we can thus obtain an upper bound of
$$ \ll_\eps \frac{1}{N^{k/2} \log^k N} \frac{1}{\log^{k \log(k-1) + o(1)} N}$$
on the probability distribution function of the tuple $(\mathbf{d}_{i,j}, \mathbf{p}_{i,j})_{1 \leq i < j \leq k}$ restricted to the event $E$ (which now forces both \eqref{dpi} and \eqref{om-d}).  By the previous arguments, we see that on an event in $E$ with probability $\gg_\eps \frac{1}{\log^k N}$, the tuple \eqref{ni-def} consists of distinct elements in $A$ multiplying to a square.  On the other hand, since $\prod_{1 \leq i < j \leq k} \mathbf{d}_{i,j}$ is squarefree on $E$, one can check from the fundamental theorem of arithmetic that any tuple $(n_1,\dots,n_k)$ can arise from at most $O(1)$ many such tuples $(d_{i,j}, p_{i,j})_{1 \leq i < j \leq k}$.  The claim follows.  It is plausible that such a bound is in fact optimal up to the $o(1)$ error in the exponent, and could also be obtained by a more direct double counting argument based on estimating the number of tuples $(n_1,\dots,n_k)$ in $\{1,\dots,N\}$ that multiply to a square and each have the ``typical'' number $(1+o(1)) \log\log N$ of prime factors, but we have not attempted to perform this calculation in detail.
\end{remark}

\begin{remark}[Higher powers]\label{higher-powers}  The above arguments can be generalized to the setting considered in the recent paper \cite{fjkps}, in which the squares are replaced by $m^{\mathrm{th}}$ powers for some fixed $m \geq 2$.  Indeed, if we let $F_{k,m}(N)$ denote the cardinality of the largest subset of $\{1,\dots,N\}$ that does not contain $k$ distinct elements multiplying to an $m^{\mathrm{th}}$ power, and set
$$ c_{k,m}^- \coloneqq \liminf_{N \to \infty} \frac{F_{k,m}(N)}{N}$$
then one can modify the proof of Theorem \ref{main-thm} to show that $c_{k,m}^- > 0$ whenever $m \geq 2$ and $k \geq m+2$ (note that the latter condition is also equivalent to $\binom{k}{m} > k$).  We sketch the relevant modifications to the proof as follows.  Firstly, we allow all constants in the asymptotic notation, as well as the parameter $\eps$ to depend on $m$ as well as $k$.  The main task is to establish a generalization of Proposition \ref{prob-prop} in which the squares in conclusion (i) are replaced by $m^{\mathrm{th}}$ powers.  Whereas in the previous argument we introduced random variables $\mathbf{d}_{i,j}$ and $\mathbf{p}_{i,j}$ for the $\binom{k}{2}$ different pairs $1 \leq i < j \leq k$, we now introduce random variables $\mathbf{d}_{i_1,\dots,i_m}$ and $\mathbf{p}_{i_1,\dots,i_m}$ for the $\binom{k}{m}$ different pairs $1 \leq i_1 < \dots < i_m \leq k$, again as independent copies of $\mathbf{d}$ and $\mathbf{p}$, but with the quantity $k-1$ appearing in the definition of $\mathbf{d}$ replaced with $\binom{k-1}{m-1}$.  A similar change is then made to the estimates \eqref{ab}, \eqref{d-dist}, \eqref{p-dist}.  We extend the notation $\mathbf{d}_{i_1,\dots,i_m}$ and $\mathbf{p}_{i_1,\dots,i_m}$ by permutations, thus for instance $\mathbf{d}_{i_{\sigma(1)},\dots,i_{\sigma(m)}} \coloneqq \mathbf{d}_{i_1,\dots,i_m}$ for any permutation $\sigma \colon \{1,\dots,m\} \to \{1,\dots,m\}$, and then define
$$ \mathbf{n}_i \coloneqq \prod_{1 \leq i_1 < \dots < i_m \leq k; i \in \{i_1,\dots,i_m\}} \mathbf{d}_{i_1,\dots,i_m} \mathbf{p}_{i_1,\dots,i_m}$$
for $i=1,\dots,k$.  The rest of the argument proceeds with the obvious notational changes, for instance the target magnitudes $N_i$ appearing before \eqref{nie} are now defined as
$$ N_i \coloneqq N / \prod_{1 \leq i_1 < \dots < i_{m-1} \leq k; i_1,\dots,i_{m-1} \neq i} d_{i,i_1,\dots,i_{m-1}}$$
and the linear forms in \eqref{u-form} are now defined as
$$ (u_{i_1,\dots,i_m})_{1 \leq i_1 < \dots < i_m \leq k} \mapsto \sum_{1 \leq i_1 < \dots < i_{m-1} \leq k; i_1,\dots,i_{m-1} \neq i} u_{i,i_1,\dots,i_{m-1}}$$
using the convention that the $d_{i_1,\dots,i_m}$ and $u_{i_1,\dots,i_m}$ are invariant under permutations of the indices.  The condition $d \geq k+2$ ensures that these linear forms remain linearly independent.  Similarly, after fixing $d_{1,l_1,\dots,l_{m-1}}$ and $p_{1,l_1,\dots,l_{m-1}}$ for $2 \leq l_1 < \dots < l_{m-1} \leq k$, the target magnitudes $\tilde N_l$ are now defined as
$$ \tilde N_l \coloneqq \frac{N}{\prod_{2 \leq l_1 < \dots < l_{m-1} \leq k: l \in \{l_1,\dots,l_{m-1}\}} p_{1,l_1,\dots,l_{m-1}} \prod_{1 \leq j_1 < \dots < j_m \leq k: l \in \{j_1,\dots,j_m\}} d_{j_1,\dots,j_m}},$$
and the linear forms replaced by
$$  (\tilde u_{i_1,\dots,i_m})_{2 \leq i_1 < \dots < i_m \leq k} \mapsto \sum_{2 \leq i_1 < \dots < i_{m-1} \leq k; i_1,\dots,i_{m-1} \neq i} \tilde u_{l,i_1,\dots,i_{m-1}},$$
which will be linearly independent for $d \geq k+2$ by similar arguments to before.  The quantities $\binom{k}{2}$ and $\binom{k-1}{2}$ that appear in these arguments will be replaced by $\binom{k}{m}$ and $\binom{k-1}{m}$ respectively, but otherwise there are no substantive changes to the arguments; we leave the details to the interested reader.
\end{remark}

\end{document}